\documentclass[a4paper,11pt]{article} 

\usepackage{amsmath,amssymb,amsthm}

\newtheorem{theorem}{Theorem}[section]

\theoremstyle{definition}
\newtheorem{definition}[theorem]{Definition}
\newtheorem{proposition}[theorem]{Proposition}

\theoremstyle{remark}
\newtheorem{remark}[theorem]{Remark}

\newcommand{\CC}{{\mathbb{C}}}
\newcommand{\PP}{{\mathbb{P}}}
\newcommand{\QQ}{{\mathbb{Q}}}

\newcommand{\Pic}{\operatorname{Pic}}

\newcommand{\smallet}{\mbox{\tiny\'et}}

\DeclareMathOperator{\tr}{trace}
\DeclareMathOperator{\Frob}{Frob}

\begin{document}

\title
{On the modularity of Calabi-Yau threefolds containing elliptic 
ruled surfaces}
\author{K.~Hulek and H. Verrill}
\date{}

\maketitle

\begin{abstract}
In this paper we prove that (not necessarily rigid) Calabi-Yau threefolds defined over $\QQ$
which contain sufficiently
many elliptic ruled surfaces are modular (under certain mild restrictions on the primes of bad
reduction). This is an application of the results of Dieulefait
and Manoharmayum \cite{DM}, \cite{D} who showed modularity for rigid Calabi-Yau threefolds.  
\end{abstract}

\section{Introduction}
In the last two years different authors found several examples of non-rigid
Calabi-Yau varieties defined over $\QQ$ for which modularity could be proved
(see \cite{LY}, \cite{HV} and \cite{Schu}). The examples of Livn\'e and Yui
arise by taking a quotient of a product of a $K3$ surface with an elliptic curve, whereas
the examples constructed by Sch\"utt and the authors have the structure of a fibre product
of elliptic surfaces. In each case the motive splits into low dimensional pieces which can be
explained geometrically.

In this brief note we want to give a criterion which implies modularity if the Calabi-Yau variety
contains sufficiently many elliptic ruled surfaces. This applies to the examples given in
\cite{HV} and \cite{Schu} and can also be used to explain the part of the motive in the examples
in \cite{LY} which does not come from the transcendental lattice of the $K3$ surface.
The main ingredient of this paper is an application of a modularity result proved by Dieulefait
and Manoharmayum in their paper \cite{DM} (and its extension by Dieulefait \cite{D})
where they establish modularity criteria for rigid Calabi-Yau threefolds.

Finally we shall show how to use our criterion to construct numereous examples of modular
(birational) Calabi-Yau varieties. In some 
of these examples we have strong numerical evidence
as to the modular forms which appear. The levels 
include 10,12,16,17,21,27,28,32,35,48,55,60,64,68,90,480 and
the corresponding varieties are listed in
Tables~\ref{table:rigid_examples},\ref{table:examples_in_caseB}, 
\ref{table:modular_examples_with_h12=1} and 
\ref{table:mostly_non_rigid_examples}.

\noindent
{\bf Acknowledgements.} We are very grateful to Professors Dieulefait and Manoharmayum for
the exchange of some e-mails in which they explained some of their results and details about
$l$-adic Galois representations to us.

We would like to thank the referee for computing the levels of the weight 
$4$ modular forms appearing in Tables~\ref{table:examples_in_caseB}, 
and \ref{table:modular_examples_with_h12=1}, and the level
$480$ example in Table~\ref{table:rigid_examples}.  
We also thank Matthias Schuett for reading through our manuscript and
making several suggestions for improvements as well as for computational help with one example.

\section{A modularity criterion}

Let $X$ be a smooth projective variety with $h^{30}(X)=h^{03}(X)=1$. (In our applications this will
be a Calabi-Yau variety or a variety birational to a nodal Calabi-Yau variety.)
A (birational) elliptic ruled surface $Y$ is a smooth projective surface which is birational
to a $\PP^1$-bundle over an elliptic curve. In this case $h^{21}(Y)=h^{12}(Y)=1$ and
$h^{30}(Y)=h^{03}(Y)=0$.

\begin{definition}
Let $X$ be smooth projective and define $b :=h^{12}(X)$. Let $Y_j, j= 1,\ldots, b$
be birational elliptic ruled surfaces contained in $X$. We say that these surfaces
{\em span $H^{21}(X) \oplus H^{12}(X)$} if the natural map
$$
 H^3(X,\CC) \to \oplus_{j=1}^b H^3(Y_j,\CC)
$$
is surjective.
\end{definition}

\begin{remark}
If the surfaces $Y_j$ are as in the previous definition, then we have an exact sequence
$$
0 \to H^{30}(X) \oplus H^{03}(X) \to H^3(X,\CC) \to \oplus_{j=1}^b H^3(Y_j,\CC) \to 0.
$$
This follows since the restriction of a $(3,0)$ or a $(0,3)$ form to an elliptic ruled
surface must vanish and together with the assumption that the surfaces 
$Y_j$ span $H^{21}(X) \oplus H^{12}(X)$ this gives the claim. In particular we obtain an 
isomorphism $H^{21}(X) \oplus H^{12}(X) \cong \oplus_{j=1}^b H^3(Y_j,\CC) = 
\oplus_{j=1}^b (H^{21}(Y_j) \oplus H^{12}(Y_j))$.
\end{remark}

\begin{remark}
If $X$ is a smooth Calabi-Yau threefold and $Y$ is a smooth elliptic ruled surface contained in $X$,
then a straightforward argument from intersection theory (cf. the proof
of Proposition \ref{prop:defect}) using the adjunction formula shows that the
map $H^3(X,\CC) \to H^3(Y,\CC)$ is surjective. If $h^{12}(X)=1$ then this shows that
we are automatically in
the situation described above. 
In fact this observation is often useful when one wants to prove
modularity of Calabi-Yau threefolds with $h^{12}(X)=1$.
\end{remark}

We now assume that we are in an arithmetic situation, i.e. that the variety $X$ as well as
the subvarieties $Y_j$ are defined over $\QQ$. When we speak about a birational elliptic surface $Y$
defined over $\QQ$ we mean that there is a birational map defined over $\QQ$ from
$Y$ to $E \times \PP^1$.
In what follows the assumption that all $Y_j$ are defined over $\QQ$ can often be weakened
to the assumption that the union of the surfaces $Y_j$ is defined over $\QQ$, but for the
sake of simplicity we shall assume that all the $Y_j$ themselves are defined over $\QQ$. Then the map
$H^3(X) \to \oplus_{j=1}^b H^3(Y_j)$ can be interpreted in the $l$-adic cohomology and gives
rise to a morphism of $\operatorname{Gal}(\overline{\mathbb{Q}}/\mathbb{Q})$-modules
$H^3_{\smallet}(X,\QQ_l) \to
\oplus H^3_{\smallet}(Y_j,\QQ_l)$. (Strictly speaking we have a 
representation of the Galois group
$\operatorname{Gal}(\overline{\mathbb{Q}}/\mathbb{Q})$ on the dual vector spaces, but we shall, by misuse
of language, refer to the cohomology groups as Galois representations.)
If the surfaces $Y_j$
span $H^{21}(X) \oplus H^{12}(X)$
this gives rise to an exact sequence
of  $\operatorname{Gal}(\overline{\mathbb{Q}}/\mathbb{Q})$-modules
\begin{equation}
\label{eq:ex_sequence}
  0 \to U \to H^3_{\smallet}(X,\QQ_l) \to
  \oplus H^3_{\smallet}(Y_j,\QQ_l) \to 0
\end{equation}
where $U$ has dimension $2$.

We shall want to make use of a criterion proved by Dieulefait and Manoharmayum \cite{DM}
and further extended by Dieulefait \cite{D}
which was developed in order to prove the modularity of rigid Calabi-Yau varieties.

\begin{proposition}
\label{prop:DM}
Let $l$ be a prime and let
$$
\rho:\operatorname{Gal}(\overline{\mathbb{Q}}/\mathbb{Q})\rightarrow
\operatorname{GL}_2 (\overline{\mathbb{Q}}_l)
$$
be an irreducible continuous Galois representation, unramified outside finitely many primes,
and having determinant $\varepsilon^3_l$ where $\varepsilon_l$ is the cyclotomic character.
Furthermore assume the following:
\begin{itemize}
\item[{\rm(1)}]$\rho$ is crystalline at $l$ with Hodge-Tate weights $\{0,3\}$.
\item[{\rm(2)}]The residual representation
$$
\overline{\rho}:\operatorname {Gal}(\overline{\mathbb{Q}}/\mathbb{Q})
\rightarrow \operatorname{GL}_2 (\overline{\mathbb{F}}_l)
$$
takes values in $\operatorname{GL}_2(\mathbb{F}_l)$.
\end{itemize}
Moreover assume that one of the following
assumptions holds:
\begin{itemize}
\item[{\rm(DM1)}]$l=7$ and $\rho$ is unramified at $3$.
\item[{\rm(DM2)}]$l=5$ and there is a prime $p\equiv \pm 2
(\operatorname{mod}5)$ where $\rho$ is unramified and the trace
of Frobenius $t_3(p)$  at $p$ is not divisible by $5$, where
we define $t_3(p):=\tr(\rho(\Frob_p))$.
\item[{\rm(D)}] $l > 3$, $\rho$ is unramified at $3$ and $t_3(3)$ is not divisible by $3$.
\end{itemize}
Then $\rho$ is modular.
\end{proposition}
\begin{proof}
Dieulefait and Manoharmayum proved that either of the conditions (DM1) or (DM2) are
sufficient for modularity (see \cite[Remark 3.2]{DM}). The fact that condition (D) is
sufficient was shown by Dieulefait (\cite[Theorem 1.1]{D}). \hfill
\end{proof}

We now want to apply this to our situation.

\begin{proposition}
\label{prop:modular}
Let $X$ be a smooth projective threefold defined over $\QQ$ with $h^{30}(X)=h^{03}(X)=1$ and assume that
$X$ contains birational ruled surfaces $Y_j, j=1, \ldots ,b$ which are
defined over $\QQ$ and which span $H^{21}(X) \oplus H^{12}(X)$. Let $\rho$ be the $2$-dimensional
Galois representation given by the kernel $U$ from the exact sequence
$$
 0 \to U \to H^3_{\smallet}(X,\QQ_l) \to
  \oplus H^3_{\smallet}(Y_j,\QQ_l) \to 0.
$$
Assume that one of the
following conditions holds:
\begin{itemize}
\item[{\rm(i)}] $X$ has good reduction at $3$ and $7$ or
\item[{\rm(ii)}] $X$ has good reduction at $5$ and some at prime
$p\equiv \pm 2 (\operatorname{mod}5)$ and the trace of $\operatorname{Frob}_p$ on $U$ is not
divisible by $5$ or
\item[{\rm(iii)}] $X$ has good reduction at $3$ and the trace of $\operatorname{Frob}_3$ on $U$ is not
divisible by $3$.
\end{itemize}
Then $X$ is modular. More precisely
\begin{equation}
\label{eq_L_series}
L(X,s) \circeq L(f_4,s) \prod_{j=1}^b L(g_2^j,s-1)
\end{equation}
where $f_4$ is a weight $4$ form and the $g_2^j$ are the weight $2$ forms associated to the base
curves
$E_j$ of the birational ruled surfaces $Y_j$. (As usual $\circeq$ denotes equality of the $L$-series
upto possibly finitely many factors.)
\end{proposition}
\begin{proof}
Since the surfaces $Y_j$ are birational over $\QQ$ to a product $E_j \times \PP^1$ we have
$$
L(Y_j,s) \circeq L(g_2^j,s-1)
$$
where $g_2^j$ is the weight $2$ form associated to the elliptic curve $E_j$.
After passing to the semi-simplification we have an isomorphism
$$
H^3_{\smallet}(X,\QQ_l) \cong U \oplus
\bigoplus_{j=1}^b H^3_{\smallet}(Y_j,\QQ_l).
$$
Our claim follows if we can show that we can apply
the results of Dieulefait and Manoharmayum, resp. Dieulefait to deduce the
modularity of the $2$-dimensional
$\operatorname {Gal}(\overline{\mathbb{Q}}/\mathbb{Q})$ representation $U$.
Clearly $U$ is
a continuous representation which is unramified outside finitely many
primes, namely the primes of bad reduction. It follows from the Weil
conjectures as proved by Deligne that $U$ is irreducible
\cite[pp. 1,2]{DM}. Moreover if $U$ is unramified at $l$ then its determinant is $\varepsilon^3_l$.
Condition $\rm{(2)}$ of \ref{prop:DM}
follows since $U$ is a subspace of $H^3_{\smallet}(X,\QQ_l)$
for which $\rm{(2)}$ holds. Similarly we can deduce that $U$ is crystalline and has
Hodge-Tate weights $\{0,3\}$. The first follows e.g. from the fact that
$H^3_{\smallet}(X,\QQ_l)$ is crystalline for $l > 3$
(see \cite[Theorem 2.2 A (ii)]{FM}) and the fact that this property is preserved for submodules
and quotients (see \cite[Proposition 1.7]{FM}). The statement about the Hodge-Tate weight
follows from standard comparison theorems (see \cite[Section 3.1]{FM}). The conditions
$\rm{(DM1)}$, $\rm{(DM2)}$ or $\rm{(D)}$ follow immediately from the hypothesis of the theorem.
\hfill
\end{proof}

\begin{remark}
The weight $4$ form which is associated to the $2$-dimensional re\-pre\-sen\-ta\-tion $U$
can in many cases be determined explicitly by the method of Faltings, Livn\'e and Serre.
\end{remark}

\section{Fibre products}

We shall now investigate a certain class of Calabi-Yau threefolds 
which often (but not always)
contain many elliptic ruled surfaces in the sense of our
discussion in the previous section.

Let $r:Y\rightarrow\mathbb{P}^1$ and $r':Y'\rightarrow\mathbb{P}^1$ be
semi-stable rational elliptic fibrations with a section. Then the fibre
product
$$
W=Y\times_{\mathbb{P}^1}Y'
$$
is a nodal Calabi-Yau variety. We denote by $\hat{W}$ a small
and by $\tilde{W}$ a big resolution of $W$. Note that small resolutions
need not be projective, in fact in many cases there exists no projective
small resolution. However,
$\hat{W}$ is a Moishezon manifold and hence we
still have a Hodge decomposition and thus can
speak about the Hodge numbers $h^{pq}(\hat{W})$, cf. \cite{Ue}.
Let $S$ and $S'$ be the sets of critical values of $r$ and $r'$.
We define $S''=S\cap S'$ and $\tilde{S}=S\cup S'\setminus
S''$. Then all fibres over points in $\tilde{S}$ are of type
$I_n\times I_0$ or $I_0\times I_n$. For any $t\in \mathbb{P}^1$
we denote by $b(t)$, resp. $b'(t)$ the number of components
of $r^{-1}(t)$, resp. $r^{-1}(t')$. Then it is well known
\cite[p. 190]{Scho} that
$$
h^{12}(\tilde{W})=h^{12}(\hat{W})= 1 +\operatorname{rank}\Pic
(W_{\eta})-\# S'' +\sum_{t\in \tilde{S}}(b(t)b'(t)-1).
$$
Here $\Pic (W_{\eta})$ is the Picard group of the general fibre.  Its
rank may be computed using the formula
$$\operatorname{rank} \Pic (W_{\eta})= \operatorname{rank} \Pic (Y_{\eta})
+
\operatorname{rank} \Pic (Y'_{\eta})+d,$$ where $d=1$ if $Y_{\eta}$ and
$Y_{\eta}'$ are isogenous and $0$ otherwise.
By \cite[p. 190]{Scho} this can be rewritten as
$$
\operatorname{rank} \Pic (W_{\eta})= d+18-\sum_{t\in S}b(t)-\sum_{t'\in S'}
b'(t')+\#S+\#S'.
$$
We define the {\it fibre defect}
$$
\delta=\delta(r,r'):= 1+\operatorname{rank} \Pic (W_{\eta})-\# S''
$$
which can be rewritten as
\begin{equation}
\delta= d + 19 -\sum_{t\in S}b(t)-\sum_{t'\in S'} b'(t')+\#S+\#S'-\#S''.
\label{eqn:formularforb}
\end{equation}
When we discuss examples in \S~\ref{subsec:newexamples}, 
we will restrict consideration
to the case where $S$ and $S'$ are rational elliptic fibrations. 
Since by Noether's theorem the Euler number of these surfaces
is $12$ it follows that
$\sum_{t\in S} b(t)=\sum_{t'\in S'} b'(t') =12$, 
and
so in this situation (\ref{eqn:formularforb}) becomes
\begin{equation}
\delta= d -5 +\#S+\#S'-\#S''.
\label{eqn:fibredefect}
\end{equation}
Using its definition, the fibre defect can be computed explicitly for all
pairs $(r, r')$ of elliptic fibrations. Note that this number is not an invariant
of the variety $W$ but depends on its representation as a fibre product (as we 
shall see in examples). We  call $\delta$ the {\it fibre defect} as this measures
the size of $H^{12}(\hat{W}) \oplus H^{21}(\hat{W})$ which does not come from fibres of
type $I_n \times I_0$.

We shall now consider the fibres over the points in $\tilde{S}$. As we
have already remarked these are of type $I_n\times I_0$ or
$I_0\times I_n$, i.e. a cycle of elliptic ruled surfaces whose components
are isomorphic to $\mathbb{P}^1_i \times E_s$ where $i=1,\dots, b(t)$,
resp. $i=1,\dots, b'(t)$ and $t\in \tilde{S}$. We shall write
$\gamma (t)$ for $b(t)$ or $b'(t)$ in order to simplify the notation.
Let $\alpha_t, \beta_t \in H_1 (E_t, \mathbb{Z})$ be generators with
$\alpha^2_t=\beta^2_t=0$ and $\alpha_t.\beta_t=1$. Then $A^i_t=
\mathbb{P}^1_i\times \alpha_t$ and $B_t^i=
\mathbb{P}^1_i\times \beta_t$ are generators
of $H_3(\mathbb{P}^1_i\times E_t)$. Finally let $\gamma_{i,t}:\mathbb{P}^1_i
\times E_t\rightarrow \hat{W}$, resp. $\gamma_{i,t}:\mathbb{P}^1_i
\times E_t\rightarrow \tilde{W}$ be the inclusion.
\begin{proposition}
\label{prop:defect}
There is an exact sequence
$$
0\rightarrow U\rightarrow H^3(\hat{W}, \mathbb{C})\overset{\oplus \gamma_{i,t}^{\ast}}
{\longrightarrow}{\bigoplus\limits_{t\in \tilde{S}}}\bigoplus\limits_{i=1}^{\gamma(t)-1}
H^3(\mathbb{P}_i^1
\times E_t, \mathbb{C})\rightarrow 0
$$
where $\operatorname{dim}_{\mathbb{C}}U=2\delta+2$. 
A similar exact sequence exists for
$\tilde{W}$ instead of $\hat{W}$.
\end{proposition}
\begin{proof}
We have to compute the rank of the intersection matrix of the cycles
$A^i_t$ and $B^i_t$. Clearly cycles in different fibres are
independent, hence it suffices to look at one fibre of type $I_n \times
I_0$. We assume that the enumeration of the components is such that
$\mathbb{P}^1_i\times E_t$ intersects $\mathbb{P}^1_j\times E_t$ if and only if
$|i-j|\leq 1$ where the index $i$ has to be read cyclically, i.e. we can
assume $i\in \mathbb{Z}/\gamma (t) \mathbb{Z}$. Then it follows immediately from
the geometry that
$$
A^i_t.A^j_t= B^i_t.B^j_t= 0 \; \mbox{ for all }i,j
$$
and
$$
A^i_t.B^j_t= -B^j_t.A^i_t=
\left\lbrace \begin{array}{lll}
0 &\mbox{ if }  & |i-j|\geq 2 \\
 1 &\mbox{ if}  & |i-j|=1.
\end{array}\right.
$$
It remains to determine $A^i_t.B^i_t$. Recall that the normal bundle
of a component $\mathbb{P}^1_i \times E_t$ in $\hat{W}$, resp. $\tilde{W}$
has degree $-2$ on the fibres $\mathbb{P}^1_i \times \left\lbrace  pt\right\rbrace$.
It then follows from \cite[19.2.2]{F} that
$$
A^i_t.B^i_t = -B^i_t.A^i_t = -2.
$$
Altogether the intersection matrix of the cycles
$A^i_t, B^i_t; i=1,\dots, \gamma (t)- 1$ is given
by
$$
\left( \begin{array}{llllll}
 \begin{pmatrix}
0 &-2  \\
2 &0
\end{pmatrix} &
 \begin{pmatrix}
 0&1  \\
-1 &0
\end{pmatrix}  &  &  &  &  \\
 \begin{pmatrix}
0 &1  \\
-1 & 0
\end{pmatrix} &
\begin{pmatrix}
0 &  -2\\
2 & 0
\end{pmatrix}
\begin{pmatrix}
 0&1  \\
 -1&0
\end{pmatrix} &  &  &  &  \\
 &  \begin{pmatrix}
 0&1  \\
 -1&0
\end{pmatrix}
\begin{pmatrix}
 0& -2 \\
2 & 0
\end{pmatrix}  &&  &  &  \\
 &  &  &  \ddots&  &  \\
 &  &  &  &
\begin{pmatrix}
 0&-2  \\
 2&0
\end{pmatrix}   &
\begin{pmatrix}
0 &1  \\
 -1&0
\end{pmatrix} \\
 &  &  &  &
\begin{pmatrix}
0 &1  \\
 -1&0
\end{pmatrix}  &
\begin{pmatrix}
 0&-2  \\
2 &0
\end{pmatrix}
\end{array}
\right)
.$$
It is obvious that this matrix has rank $2\gamma (t)-2$ and this proves the proposition.
\hfill
\end{proof}

\begin{remark}
It is actually not necessary that we start with semi-stable fibrations provided $\hat{W}$ is
a Calabi-Yau variety and we
obtain sufficiently many contributions from the fibres of type $I_n \times I_0$ .
\end{remark}

We now want to consider the modularity of such fibre products.

\begin{theorem}
\label{modularity}
Assume that $r:Y\rightarrow \mathbb{P}^1$ and $r':Y'\rightarrow \mathbb{P}^1$
are semi-stable, rational elliptic fibrations with a section, which are defined
over $\mathbb{Q}$. Let $\tilde{W}$ be a (big) resolution of
$W=Y \times_{\mathbb{P}^1}Y'$. Moreover assume that
\begin{itemize}
\item[{\rm(1)}] the fibre defect $\delta=\delta(r,r')=0$,
\item[{\rm(2)}] if $t \in \tilde{S}$ and the fibre $W_t$ is isomorphic to
$E_t\times I_n$ for some elliptic curve $E_t$ and an $n$-cycle $I_n,
n\geq 2$, then all components of this fibre are defined over $\mathbb{Q}$.
\end{itemize}
Then there is an exact sequence of $\operatorname {Gal}(\overline{\mathbb{Q}}/\mathbb{Q})$-modules
$$
  0 \to U \to H^3_{\smallet}(\tilde{W},\QQ_l) \to
  {\bigoplus\limits_{t\in \tilde{S}}}\bigoplus\limits_{i=1}^{\gamma(t)-1}
  H^3_{\smallet}(\mathbb{P}_i^1 \times E_t)\to 0
$$
with $2$-dimensional kernel $U$.

Moreover assume that one of the following assumptions holds:
\begin{itemize}
\item[{\rm(3)}]$\tilde{W}$ has good reduction at $3$ and $7$ or
\item[{\rm(4)}]$\tilde{W}$ has good reduction at $5$ and some prime
$p\equiv \pm 2 (\operatorname{mod}5)$ with the trace of $\operatorname{Frob}_p$ on the
kernel $U$ not divisible by $5$ or
\item[{\rm(5)}]$\tilde{W}$ has good reduction at $3$ and the trace of $\operatorname{Frob}_3$ on $U$
not divisible by $3$.
\end{itemize}
Then $\tilde{W}$ is modular and the same is true for a small
projective resolution defined over $\mathbb{Q}$ if it exists. Moreover
$$
L (H^3_{\smallet}(\tilde{W}),s) \circeq L(f,s) \prod\limits_{t\in\tilde{S}}
(g_{E_{t}}, s-1)^{\gamma (t)-1}
$$
where $f$ is a modular form of weight $4$ and $g_{E_t}$ is the
weight $2$ form associated to $E_t$.
\end{theorem}
\begin{proof}
By Proposition \ref{prop:defect} the condition $\delta=0$ assures that the
components of the fibres $E_t \times \PP^1$ over the
points in $\tilde S$ span $H^{21}(\tilde W) \oplus H^{12}(\tilde W)$ and that the kernel
$U$ has, therefore, dimension $2$. The result is now an
immediate application of Proposition \ref{prop:modular}.
\hfill
\end{proof}

\section{Examples}

In this section we want to discuss various examples where the method described above can
be applied.  We use
(\ref{eqn:fibredefect}) to find the fibre defect $\delta$ in these
examples.

\subsection{Relation with previous examples}

In \cite{HV} we considered the following examples. Let
$\mathcal E_{(\alpha,\beta,\gamma)}$ be the elliptic fibration which is the resolution
of the surface $\mathcal E'_{(\alpha,\beta,\gamma)}\subset\PP^2\times\PP^1$
given by
\begin{eqnarray}
\label{eqn:Elliptic_abct_family}
\mathcal E'_{(\alpha,\beta,\gamma)}:&&(x+y+z)(\alpha xy + \beta yz + \gamma zx)t_0=t_1xyz,
\end{eqnarray}
where $(x:y:z)\in\PP^2, (t_0:t_1)\in\PP^1$.
There is a projection
$p:\mathcal E_{(\alpha,\beta,\gamma)}\rightarrow\PP^1$, with
fibre $\mathcal E_{(\alpha,\beta,\gamma),t}:=p^{-1}(1:t)$ for generic $t$.
We write $\mathcal E_{\alpha}:=\mathcal E_{(1,1,\alpha)}$
and $\mathcal E_{\alpha,t}:=\mathcal E_{(1,1,\alpha),t}$.
The singular fibres of these families can be read off from \cite[Section 4, Table 2]{HV}.

Taking fibre products and resolving singularities one obtains smooth varieties
$$
X_{\bf a} =  X_{({\alpha}_1:{\alpha}_2:{\alpha}_3:{\alpha}_4:{\alpha}_5:{\alpha}_6)} \sim_{bir}
{\mathcal E}_{({\alpha}_1:{\alpha}_2:{\alpha}_3)} \times_{{\mathbb P}^1}
{\mathcal E}_{({\alpha}_4:{\alpha}_5:{\alpha}_6)}
$$
which are resolutions of nodal Calabi-Yau threefolds. The examples discussed
in \cite{HV} were  $X_{25}=X_{(1:1:1:1:1:25)}$, $X_{(1:1:1:9:9:9)}$
and $X_{(1:1:1:4:4:9)}$ where $h^{12}=4,2$ and $1$ respectively.
In the last case
$$
X_{(1:1:4:4:4:16)} \sim_{bir} {\mathcal E}_{(1:1:16)} \times_{{\mathbb P}^1} {\mathcal E}_{(4:4:4)}
$$
and we have  $\#S=  \#S' =4$ and $ \#S''=3$ and hence, since these families are not isogeneous, we
have  $\delta=0$. In this case $h^{12}=1$ and it is
sufficient to consider the
$I_2 \times I_0$-fibre over the point $16$. We can then derive modularity
from Theorem \ref{modularity}.

In the second example we can either write
$$
X_{(1:1:1:9:9:9)} \sim_{bir} {\mathcal E}_{(1:1:1)} \times_{{\mathbb P}^1} {\mathcal E}_{(9:9:9)}
$$
or
$$
X_{(1:1:1:9:9:9)} \sim_{bir} {\mathcal E}_{(1:1:9)} \times_{{\mathbb P}^1} {\mathcal E}_{(1:9:9)}.
$$
In the first case we have again $\#S=  \#S' =4$ and $ \#S''=3$ and hence  $\delta=0$. We have
a fibre of type $I_3 \times I_0$ over the point $1$ and since $h^{12}=2$ this spans all of
$H^{21} \oplus H^{12}$ and our method applies. In the second description we have  $\#S=  \#S' =5$
and $ \#S''=4$ and this gives $\delta=1$. This is also an example which shows that $\delta$ is not a 
birational invariant but depends on the chosen fibration.

Finally in the first example
$$
X_{25} \sim_{bir} {\mathcal E}_{(1:1:1)} \times_{{\mathbb P}^1} {\mathcal E}_{(1:1:25)}
$$
and in this case  $\#S=4$,   $\#S'=5$, $\#S''=3$ and thus $\delta=1$. Hence we cannot invoke
Theorem \ref{modularity}. However, we showed in \cite{HV} that one can still find
sufficiently many elliptic ruled surfaces in $X_{25}$ in order to be able to apply
Proposition \ref{prop:modular}.

The weight $4$ forms corresponding to $H^{30} \oplus H^{03}$ were computed in \cite{HV}.

A further class of examples where this method works was discussed by Sch\"utt \cite{Schu}.
In his examples  $\#S=5$ or $4$ and $\#S'=4$ with $\#S''=4$ and $3$ respectively. In all
cases this gives $\delta=0$ and hence one is in a situation where Theorem \ref{modularity}
applies.

\subsection{New examples}
\label{subsec:newexamples}

In this section we give a large number of examples of (birational)
Calabi-Yau threefolds
for which
Theorem~\ref{modularity} implies modularity.  
In many cases we shall also conjecture the modular forms which appear.
We have strong numerical evidence for these conjectures, i.e. we computed some of the
Fourier coefficients modulo $p$ for all 
primes of good reduction less than $100$.
It seems very likely that, using the technique of
Faltings, Livn\'e and Serre, these conjectures can be verified.

We use the following notation.
If $\mathcal F^t$ is an elliptic fibration with pa\-ra\-me\-ter $t$,
and if $M$ is a fractional linear transformation, then
$\mathcal F^{Mt}$ denotes the family obtained by replacing $t$ with $Mt$.

We denote by $\mathcal E(\Gamma)$ 
Beauville's modular elliptic surfaces
\cite{B}, where
$\Gamma$ is one of
$\Gamma(3)$
$\Gamma_1(4)\cap\Gamma(2)$,
$\Gamma_1(6)$,
$\Gamma_1(5)$,
$\Gamma_0(8)\cap\Gamma_1(4)$ or
$\Gamma_0(9)\cap\Gamma_1(3)$.
Note that $\mathcal E(\Gamma_1(6))=\mathcal E_{(1:1:1)}$.
The others may be given by equations as in the table below.
The values of $t$ for which the fibres are
singular are tabulated, and the last column shows the
number of components of the singular fibres, all of which are
of $I_n$ type.
$$
\begin{array}{l|l|l|l}
\Gamma &\text{equation for }
\mathcal E(\Gamma) & \text{singular }t &
\text{fibres}\\
\hline
\Gamma(3) &   x^3+y^3+z^3=3t xyz&
\infty,1,\frac{-1\pm\sqrt{-3}}{2}
&
3,3,3,3\\
\Gamma_1(4)\cap\Gamma(2)& 
x(x^2 + z^2 + 2zy)=tz(x^2-y^2)
&\infty,0,1,-1&
4,4,2,2\\
\Gamma_1(5)           & x(x-z)(y-z)=tzy(x-y)
&\infty,0,\frac{-11\pm5\sqrt{5}}{2}&
5,5,1,1\\
%
\Gamma_1(6)           & (x+y+z)(xy+yz+zx)=txyz
&\infty,0,1,9&
6,2,3,1\\
\Gamma_0(8)\cap\Gamma_1(4)& 
(x+y)(xy+z^2)=4t xyz&
\infty,0,1,-1&
8,2,1,1\\
\Gamma_0(9)\cap\Gamma_1(3)& x^2y + y^2z+z^2x=3t xyz
& \infty,1,
\frac{-1\pm\sqrt{-3}}{2}&
9,1,1,1
\end{array}
$$

If $Y$ and $Y'$ are {\bf not} isogenous then 
the condition $\delta=0$ 
in Theorem~\ref{modularity}
translates into 
$$
\#S + \#S' = \#S'' +5.
$$
This shows that we cannot have examples where at least $1$ elliptic fibration has $6$ singular
fibres. 
In order to produce examples with $\delta=0$ we have three possibilities:

$$
\renewcommand\arraystretch{1.5}
\begin{array}{l|ccc|l}
& \#S & \#S' & \#S'' &
\multicolumn{1}{|c}{\text{possible examples}}\\
\hline
\text{Case }A & 5 & 5 & 5 & 
\mathcal E_{(1:\alpha^2:\alpha^2)}\times_{\PP^1}
\mathcal E_{(1:\beta^2:\beta^2)}^{Mt}\\
\hline
\text{Case }B & 5 & 4 & 4 &
\mathcal E_{(1:\alpha^2:\alpha^2)}\times_{\PP^1}
\mathcal E_{(1:1:1)}^{Mt}\\
&&&&
\mathcal E_{(1:\alpha^2:\alpha^2)}\times_{\PP^1}
\mathcal E^{Mt}(\Gamma)\\
\hline
\text{Case }C & 4 & 4 & 3 &
\mathcal E_{(1:1:1)}\times_{\PP^1}
\mathcal E_{(1:1:1)}^{Mt}\\
&&&& 
\mathcal E(\Gamma_1)\times_{\PP^1}\mathcal E^{Mt}(\Gamma_2)
\end{array}
$$

We now construct explicit examples for each case.

\noindent
{\bf Case A:} $\#S = \#S' = \#S'' =5$ (and the two surfaces are not isogeneous).
Notice that this can never give non-rigid 
examples since we have no fibres of 
type $I_n \times I_0$. 

Given a fibration  $\mathcal E_{(\alpha,\beta,\gamma)}$ and a fractional linear 
transformation $M$, we consider
$$
\mathcal E_{(\alpha,\beta,\gamma)}^{Mt}: (x+y+z)(axy + byz +czx)=(Mt)xyz.
$$
We shall take two families 
given by
$\mathcal E_{(1,\alpha^2,\alpha^2)}$
and $\mathcal E_{(1,\gamma^2,\gamma^2)}$,
with $1,\alpha,2\alpha,0$ distinct, and 
$1,\gamma,2\gamma,0$ distinct, and analyse 
possible fibre products. Recall from \cite[Section 4, Table 2]{HV} 
that the singular fibres of $\mathcal E_{(1,\alpha^2,\alpha^2)}$  
lie over the points $\{\infty,0,1,(1 \pm 2\alpha)^2\}$, with
fibre types
$I_6,I_2,I_2,I_1,I_1$ respectively, provided $2\alpha\not=0,1,2$.
If there is a fractional linear transformation $M$ with
$M\{\infty, 0, 1, (1\pm2\alpha)^2\}= 
\{\infty, 0, 1, (1\pm2\gamma)^2\}$,
then the resolution of
$$\mathcal E_{(1,\alpha^2,\alpha^2)}\times
\mathcal E_{(1,\gamma^2,\gamma^2)}^{Mt}$$
is a rigid (birational) Calabi-Yau threefold, 
provided that $\#\{\infty,0,1,(1\pm2\alpha)^2\}=5$
and that this is not an isogeneous pairing.
There are $60$ such choices 
(this can be reduced by symmetry considerations)
for which three of the points of
$\infty, 0, 1, (1\pm2\gamma)^2$ map to $0$, $1$ and $\infty$.
This determines $M$, and the remaining two points $\beta,\delta$ then map to
$M(\beta)$, $M(\delta)$, which are given by
rational functions in $\alpha$ and $\gamma$.
The requirement that
\begin{equation}
\label{eqn:transforms}
M(\beta)=(1+ 2\alpha)^2\; \text{ and }\; M(\delta)=(1- 2\alpha)^2,
\end{equation}
implies $\alpha=(M(\beta)-M(\delta))/8$.
This expression for $\alpha$ in terms of $\gamma$ can be
substituted into (\ref{eqn:transforms}) to give polynomials
satisfied by $\gamma$, and from this we can determine whether
$\alpha^2$ and $\gamma^2\in\QQ\setminus\{0\}$.  One must
also check that $(1\pm2\gamma)^2\not=0,1,\infty$.
In this case, we get a new example of a rigid Calabi-Yau threefold.

We list all possible cases
in Table~\ref{table:rigid_examples}.  
The singular fibres of the first
family have types $I_6, I_2, I_2, I_1, I_1$, in this
order, so we only list the fibre types of the 
second family.
\begin{table}
$$
\renewcommand\arraystretch{1.5}
\begin{array}{ccc|ccccc|cc|cc|c}
&&&\multicolumn{5}{|c|}{\text{fibre types of 
$\mathcal E^{Mt}_{(1, \gamma^2,\gamma^2)}$
}}
&&&&&\text{level}\\
\alpha^2 & \gamma^2 & Mt
& \infty & 1 & 0 &
\multicolumn{2}{c|}{(1\pm 2\alpha)^2}
&\multicolumn{2}{|c}{(1\pm 2\alpha)^2}
&\multicolumn{2}{|c|}{(1\pm 2\gamma)^2}\\
\hline
\frac{-1}{8}  & \frac{-1}{8} &  1-t
& I_6 & I_2& I_2& I_1& I_1 & 
\multicolumn{2}{|c|}{\frac{1\pm\sqrt{-8}}{2}}  &
\multicolumn{2}{|c|}{\frac{1\pm\sqrt{-8}}{2}}  &32
\\
\frac{1}{2} & \frac{1}{2} & 1/t
& I_2 & I_6 & I_2 & I_1 & I_1 &
\multicolumn{2}{|c|}{3\pm2\sqrt{2}}&
\multicolumn{2}{|c|}{3\pm2\sqrt{2}}&64
\\
\frac{1}{9} & 4 & 9t
& I_6 & I_2 & I_1 & I_2 & I_1 
& 1/9 & 25/9 & 9 & 25&60
\\
\frac{9}{4} & \frac{9}{4} &  16/t
& I_2 & I_6 & I_1 & I_2 & I_1
 &16 & 4 &16 & 4 & 480
\\
\frac{1}{64} & \frac{1}{64} & \frac{25}{16} -t
& I_6 & I_1 & I_1 & I_2 & I_2 
&\frac{9}{16} & \frac{25}{16}&\frac{9}{16} & \frac{25}{16}&35
\end{array}
$$
\caption{Rigid Calabi-Yau threefolds of the form
$\mathcal E_{(1,\alpha^2,\alpha^2)}\times
\mathcal E_{(1,\gamma^2,\gamma^2)}^{Mt}$
}
\label{table:rigid_examples}
\end{table}
Note that
$\mathcal E^t_{(1,1/9,1/9)}\times\mathcal E^{9t}_{(1,4,4)}$
is the same as
$
\mathcal E_{(9,1,1)}\times\mathcal E_{(1,4,4)}
\sim_{bir} 
X_{(1:1:1:4:4:9)}$, which we considered in \cite{HV},
so this case gives only four new cases.

Computational evidence indicates that the L-series of the 
middle cohomology of
$\mathcal E^t_{(1,-1/8,-1/8)}\times\mathcal E^{1-t}_{(1,-1/8,-1/8)}$
is a weight $4$ level $32$ modular form, with 
$q$-expansion
$$
q - 8q^3 - 10q^5 - 16q^7 + 37q^9 + 40q^{11}  +\cdots.$$

Similarly, computations indicate that 
$\mathcal E^t_{(1,\frac{1}{2},\frac{1}{2})}\times\mathcal E^{1/t}_{(1,\frac{1}{2},
\frac{1}{2})}
$
is associated with a weight $4$ level $64$ modular form, with 
$q$-expansion
$$q - 22q^5 - 27q^9 + 18q^{13}+\cdots,$$
$\mathcal E^t_{(1,9/4,9/4)}\times\mathcal E^{16/t}_{(1,9/4,9/4)}$
is associated with a weight $4$ level $480$ modular form, 
$$q + 3q^3 + 5q^5 + 4q^7 + 9q^9 - 40q^{11} - 90q^{13}+\cdots,$$
and 
$\mathcal E^t_{(1,\frac{1}{64},\frac{1}{64})}
\times\mathcal E^{25/16-t}_{(1,\frac{1}{64},\frac{1}{64})}
$
is associated with
the weight $4$ level $35$ form 
$$
q + q^2 - 8q^3 - 7q^4 - 5q^5 - 8q^6 +\cdots
.
$$

\noindent
{\bf Case B:} $\#S = 5$, $\#S' = \#S'' =4$.
\nopagebreak
\newline
\noindent
{\bf i. Fibre products of $\mathcal E_{(1:a:a)}$ 
fibrations.}
We first consider products of the form
$$
{\mathcal E}_{(1:\alpha^2:\alpha^2)} 
\times_{{\mathbb P}^1} {\mathcal E}_{(1:1:1)}^{Mt}
$$
for some fractional linear transformation $M$, and 
$\alpha\not=0,1,1/2$.
Since we are more interested in non-rigid cases, we will
not investigate the rigid possibilities. Recall that the singular fibres of ${\mathcal E}_{(1:1:1)}$
lie over the points $\{\infty,0,1,9\}$ and that the singular fibres are of type $I_6,I_2,I_3,I_1$
(in this order).
For a non-rigid example, we are therefore looking for matrices $M$ with
$$
\{\infty,0,1\}\not\subset
M\{\infty,0,1,9\}\subset\{\infty,0,1,(1\pm2\alpha)^2\}.
$$
Denote the cross ratio 
by $C(a,b,c,d)=(a-c)(b-d)/(a-d)(b-c)$, and set $j(a,b,c,d)=j(C(a,b,c,d))$,
where
$$j(\lambda)=\frac{(\lambda^2-\lambda+1)^3}{(\lambda-1)^2
\lambda^2}.$$
Possible values of $\alpha^2$ are given by solutions to
$73^32^{-6}3^{-4}=j(\infty,0,1,9)=J$, for
$J=j(0,1,(1\pm2\alpha)^2)$,
$j(\infty,1,(1\pm2\alpha)^2)$ or
$j(0,\infty,(1\pm2\alpha)^2)$.
The first possibility,
$$j\left(\left(\frac{2\alpha+1}{2\alpha-1}\right)^2
\left(\frac{\alpha+1}{\alpha-1}\right)\right)
=73^32^{-6}3^{-4}$$
 has no rational solutions.
Up to sign, the second,
$$
j\left(\frac{\alpha+1}{\alpha-1}\right)=\frac{(\alpha^2+3)^3}{4
(\alpha^2-1)^2}=73^32^{-6}3^{-4}
$$
 gives $\alpha=17, 5/4, 7/9$,
and the third,
$$
j\left(\left(\frac{2\alpha+1}{2\alpha-1}\right)^2\right)
=73^32^{-6}3^{-4}
$$
 gives $\alpha=1, 1/4$ (plus irrational solutions),
but $\alpha=1$ is a degenerate case.
We obtain the examples given in Table~\ref{table:examples_in_caseB}.
As before, the first family,
$\mathcal E_{(1:\alpha^2:\alpha^2)}$,  has fibre types
$I_6, I_2, I_2, I_1, I_1$ in this order, so we only
list the fibre types for the second family,
$\mathcal E^{Mt}_{(1:1:1)}$.
In all cases $h^{12}=1$.  Note that the two level $90$ 
modular forms in this table are not twists of each other.
\begin{table}
$$
\begin{array}{cc|ccccc|cc|cl}
\alpha & Mt
&\multicolumn{5}{|c|}{\text{fibre types}}&&&\text{level}&
\multicolumn{1}{c}{\text{modular form}}\\
&& \infty & 1 & 0 &
\multicolumn{2}{c|}{(1\pm 2\alpha)^2}
&\multicolumn{2}{|c|}{(1\pm 2\alpha)^2}
&&\\
\hline
\rule{0ex}{3ex}
17 & \frac{(t-1)}{(t-33^2)}
& I_3 & I_2& -& I_6& I_1 & 
33^2&35^2
&17&q-3q^2-8q^3\cdots
\\
\rule{0ex}{3ex}
\frac{5}{4} & \frac{45}{4(t-1)}
& I_2 & I_6 & - & I_1 & I_3 &
\frac{9}{4} & \frac{49}{4}
&90&q-2q^2+4q^4\\
&&&&&&&&&\multicolumn{2}{r}{-5q^5-4q^7-8q^8}\\
&&&&&&&&&\multicolumn{2}{r}{+ 10q^{10} - 12q^{11}\cdots}
\\
\rule{0ex}{3ex}
\frac{7}{9} & \frac{(81t-25)}{56}
& I_6 & I_3 & - & I_2 & I_1 
& \frac{25}{81} & \frac{23^2}{81}&21&q-3q^2-3q^3\cdots
\\
\rule{0ex}{3ex}
\frac{1}{4} & \frac{9}{4t}
& I_2 & - & I_6 & I_1 & I_3
&\frac{1}{4} & \frac{9}{4} &90&
q-2q^2+4q^4 \\
&&&&&&&&&\multicolumn{2}{r}{+5q^5-4q^7-8q^8}\\
&&&&&&&&&\multicolumn{2}{r}{-10q^{10}+48q^{11}\cdots}
\end{array}
$$
\caption{Modular Calabi-Yau varieties of the form
${\mathcal E}_{(1:\alpha^2:\alpha^2)} 
\times_{{\mathbb P}^1} {\mathcal E}_{(1:1:1)}^{Mt}$
}
\label{table:examples_in_caseB}
\end{table}

\noindent
{\bf ii. $\mathcal E_{(1:a:a)}$ with Beauville modular fibrations.}
The same procedure may be applied to fibrations
of the form 
$$\mathcal E_{(1:\alpha^2:\alpha^2)}\times_{\PP^1}
\mathcal E(\Gamma)^{Mt}.$$
Table~\ref{table:modular_examples_with_h12=1}, which is not exhaustive,
gives a few examples,
all of which have $h^{12}=1$.
Note that the level $55$ example is the same as
Sch\"utt's $W_2$ in \cite{Schu}. There is a $3:1$ isogeny between the two 
cases in Table~\ref{table:modular_examples_with_h12=1}
where level $27$ occurs.  This isogeny is
described by Schoen \cite{schoen2}.

\begin{table}
$$
\begin{array}{l|l|l|ll}
\alpha^2 & \Gamma & M  & \text{level}&\text{modular form}\\
\hline
-3&\Gamma(3) & (t+7)/8 & 27& q -3q^2 + q^4 -15q^5 \cdots \\
9&\Gamma_1(4)\cap\Gamma(2) & (25-t)/24&12 & 
q + 3q^3 - 18q^5 + 8q^7\cdots\\
121/125 &\Gamma_1(5) & \frac{125}{88}(1-t) &55&
q + q^2 - 3q^3 - 7q^4+\cdots \\
9&\Gamma_0(8)\cap\Gamma_1(4) & 24/(25-t)&12 &
q + 3q^3 - 18q^5 + 8q^7\cdots\\
-3&\Gamma_0(9)\cap\Gamma_1(3)& (t+7)/8 & 27 & q -3q^2 + q^4 -15q^5 \cdots \\
\end{array}
$$
\caption{Modular Calabi-Yau threefolds of the form
$\mathcal E_{(1:\alpha^2:\alpha^2)}\times_{\PP^1}
\mathcal E(\Gamma)^{Mt}$
}
\label{table:modular_examples_with_h12=1}
\end{table}

\noindent
{\bf Case C:} $\#S = \#S' = 4$,  $\#S'' =3$. 
\nopagebreak
\newline
\noindent
{\bf i. Fibre products of $\mathcal E_{(1:1:1)}$ 
fibrations.}
We consider non-rigid possibilities of the form
$$
{\mathcal E}_{(1:1:1)} 
\times_{{\mathbb P}^1} {\mathcal E}_{(1:1:1)}^{Mt}. 
$$
We must take $M$ such that
$$\#(\{0,1,\infty,9\}\cup M\{0,1,\infty,9\})=3.$$
For example, there are $6$ possible ways to
map $0,1,9$ to $0,1,9$, only one of which (the identity)
maps $\infty$ to $\infty$.
The other $5$ cases are given by transformations
given by matrices
$$
M:=\begin{pmatrix}81& -81\\17 &-81\end{pmatrix},
\begin{pmatrix}81& -81\\73& -9\end{pmatrix},
\begin{pmatrix}9& -81\\73& -81\end{pmatrix},
\begin{pmatrix}9& 0\\10& -9\end{pmatrix},
\text{ or}
\begin{pmatrix}-1& 9\\7& 1\end{pmatrix}.
$$
The resulting fibrations all have $h^{12}=10$.

For each of the $10$
choices of pairings, we have $5$ maps for which
$S\not=S'$ (there are $6$ maps, 
(the number of permutations of three object), 
but this includes the identity, or a map
such as $t\mapsto 9/t$ which
leaves the cross ratio of $0,1,9,\infty$
invariant 
(the  cross-ratio is invariant under 
$a\leftrightarrow b, c\leftrightarrow d$)).
Explicitly, the maps are given by fractional 
linear transformations defined by matrices
$m_iT^j(R T)^km_l^{-1}$ where 
$1\le i\le l\le 4$, $j=0,1$, $k=0,1,2$, 
and
$T=\begin{pmatrix}-1 & 1\\ 0 & 1 \end{pmatrix}$,
$R=\begin{pmatrix}0 & 1\\ 1 & 0 \end{pmatrix}$,
$m_1=\begin{pmatrix}9 & 0\\ 1 & 8 \end{pmatrix}$,
$m_2=\begin{pmatrix}9 & 0\\ 0 & 1 \end{pmatrix}$,
$m_3=\begin{pmatrix}8 & 1\\ 0 & 1 \end{pmatrix}$,
$m_4=I$.
This gives $50$ fibrations,
(six of which
are given in 
Table~\ref{table:mostly_non_rigid_examples}),
which may not all be new, for example,
$$
\mathcal E_{(1:1:1)}\times\mathcal E^{9t}_{(1:1:1)}
\sim_{bir} 
X_{(1:1:1:9:9:9)},$$
and may not all be distinct, for example, 
$\mathcal E_{(1:1:1)}\times\mathcal E^{Mt}_{(1:1:1)}$
is the same variety as
$\mathcal E_{(1:1:1)}\times\mathcal E^{M^{-1}t}_{(1:1:1)}$.
The only rigid examples
are those with
$M\{0,1,\infty\}=\{0,1,\infty\}$, 
(two of which are the same) so we have
$4$ rigid cases, and at most $45$ non-rigid cases.
All of the rigid cases 
(given as the first $4$ examples in 
Table~\ref{table:mostly_non_rigid_examples})
have already been studied in detail by
Sch\"utt \cite{Schu2}.
Sch\"utt also pointed out to us that
the level $90$ example in 
Table~\ref{table:mostly_non_rigid_examples}
is a twist of the level $10$ example in this table.
Note that this level $90$ form is the same as the first
level $90$ example in Table~\ref{table:examples_in_caseB},
but for that case we had $h^{12}=1$. We believe that those examples in
Table~\ref{table:mostly_non_rigid_examples} which have the same modular form are
isogeneous.

\noindent{\bf ii. Beauville case.}

We can obtain examples from the 
$\mathcal E(\Gamma)$ in exactly the same way 
as for $\mathcal E_{(1:1:1)}=\mathcal E(\Gamma_1(6))$ above.
However, in general,
unless all the singular $t$ are rational,
we have fewer cases.  Also in some examples, many of the
possible $M$ give rise to isomorphic cases.
Table~\ref{table:mostly_non_rigid_examples}
gives a number of examples of modular fibrations
of the form
$$\mathcal E(\Gamma_1)\times_{\PP^1}
\mathcal E^{Mt}(\Gamma_2),$$
but is far from exhaustive.  The level means the conjectural
level of the  corresponding weight $4$
modular form, indicated by computations.
\begin{table}
$$
\begin{array}{lll|l|ll}
\Gamma_1 & \Gamma_2 & M & h^{12} &\text{level}&\text{modular form}\\
\hline
\Gamma_1(6) &\Gamma_1(6) & t/(t-1) & 0 & 21 & q - 3q^2  - 3q^3\cdots \\
\Gamma_1(6) &\Gamma_1(6) & 1/t  & 0 & 10 &q + 2q^2 - 8q^3 \cdots \\
\Gamma_1(6) &\Gamma_1(6) & 1-t  & 0 &17&q-3q^2-8q^3\cdots\\
\Gamma_1(6) &\Gamma_1(6) & (t-1)/t  & 0 &73&q + 3q^2-8q^3\cdots\\
\Gamma_1(6) &\Gamma_1(6) & 9t & 2 & 90 &q - 2q^2 + 4q^4\\
&&&&\multicolumn{2}{r}{-5q^5-4q^7-8q^8\cdots}\\
\Gamma_1(6) &\Gamma_1(6) & 9-t & 4 & 21 & q - 3q^2  - 3q^3\cdots \\
\Gamma_0(8)\cap\Gamma_1(4) &
\Gamma_0(8)\cap\Gamma_1(4) & 
1-t & 0 & 12 & q + 3q^3 - 18q^5\cdots \\
\Gamma_1(4)\cap\Gamma(2) &
\Gamma_1(4)\cap\Gamma(2) &
1-t & 2 & 12 & q + 3q^3 - 18q^5\cdots \\
\Gamma_0(8)\cap\Gamma_1(4) & 
\Gamma_0(8)\cap\Gamma_1(4) & \frac{-2}{t+1} & 1 & 48 &
q - 3q^3 - 18q^5 \cdots \\
\Gamma_0(8)\cap\Gamma_1(4) & \Gamma_1(6) & t &0&80&q - 2q^3 - 5q^5 \cdots\\
\Gamma_1(4)\cap\Gamma(2) & \Gamma_1(6) & t & 1&80&q - 2q^3 - 5q^5 \cdots\\
\Gamma_1(4)\cap\Gamma(2) & \Gamma_1(6) & 1/t & 1&80&q - 2q^3 - 5q^5 \cdots\\
\Gamma_1(4)\cap\Gamma(2) & \Gamma_1(6)& 1-t& 1&28& q - 10q^3 - 8q^5\cdots\\
\Gamma_1(4)\cap\Gamma(2) & \Gamma_1(6)& 1/(1-t)& 1&28& q - 10q^3 - 8q^5\cdots\\
\Gamma_1(4)\cap\Gamma(2) & \Gamma_1(6)& 1-1/t& 1&68& q - 2q^3 - 8q^5\cdots\\
\Gamma_0(9)\cap\Gamma_1(3) 
& 
\Gamma_0(9)\cap\Gamma_1(3) & 1/t & 16 & 27 & q - 3q^2 + q^4 \cdots\\
\end{array}
$$
\caption{Modular Calabi-Yau threefolds of the form
$\mathcal E(\Gamma_1)\times
\mathcal E(\Gamma_2)^{Mt}$
}
\label{table:mostly_non_rigid_examples}
\end{table}
Note that there are some isogenous pairs in this table, for which we expect that
the weight $4$ modular forms are the same.  However, $h^{12}$ is  not preserved by the
map corresponding to the isogeny.  Several cases
involving
$\Gamma_1(4)\cap\Gamma(2)$, with $h^{12}=1$,
could be replaced by rigid cases with
$\Gamma_0(8)\cap\Gamma_1(4)$ replacing $\Gamma_1(4)\cap\Gamma(2)$.

\bigskip

\noindent{\bf The isogenous case:}
Finally, 
If $Y$ and $Y'$ {\bf are} isogenous then 
$\#S = \#S' = \#S''$, and
the condition $\delta=0$ 
in Theorem~\ref{modularity}
translates into 
$\#S = 4.$  

All examples of this type are rigid.
This case includes all self fibre products of Beauville
examples,  given in the appendix to
\cite{Y}.  Another example is 
$$\mathcal E(\Gamma_0(8)\cap\Gamma_1(4))\times_{\PP^1}
\mathcal E^{(t-1)/(t+1)}(\Gamma_0(8)\cap\Gamma_1(4))$$
The fact that
$t\mapsto (t-1)/(t+1)$ 
is an isogeny is explained in more detail
in Example~3.5 in 
\cite{V-Moonshine_proc}, which deals with the same elliptic
pencil, but with parameter $t$ replaced by
$\frac{1}{4t}$.

We have computationally verified that
the Mellin transform of the L-series of the 
middle cohomology of this variety is
the weight $4$ level $16$
modular form, $$q + 4q^3 - 2q^5  - 24q^7 - 11q^9 + 44q^{11}\cdots.$$

\bibliographystyle{alpha}

\bigskip

\bigskip

\noindent
\parbox{2.7in}{
Klaus Hulek\\
Institut f\"ur Mathematik (C)\\
Universit\"at Hannover\\
Welfengarten 1, 30060 Hannover\\
 Germany\\
{\tt hulek@math.uni-hannover.de}
}
\parbox{2in}{
Helena Verrill\\
Department of Mathematics\\
Louisiana State University\\
Baton Rouge, LA 70803-4918\\
USA\\
{\tt verrill@math.lsu.edu}
}

\end{document}